\documentclass[12pt]{amsart}
\usepackage{amssymb, amstext, amscd, amsmath}
\makeatletter
\def\@cite#1#2{{\m@th\upshape\bfseries%
[{#1\if@tempswa{\m@th\upshape\mdseries, #2}\fi}]}}
\makeatother
\renewcommand{\labelenumi}{(\roman{enumi}) }
\theoremstyle{plain}
\newtheorem{thm}{Theorem}[section]

\newtheorem{lem}[thm]{Lemma}
\theoremstyle{definition}

\newtheorem{rems}[thm]{Remarks}
\newtheorem{defn}[thm]{Definition}

\newtheorem{ques}[thm]{Question}

\newcommand{\Prf}{\noindent\textbf{Proof.\ }}
\newcommand{\bx}{\hfill$\blacksquare$\medbreak}

\newcommand{\BH}{{\B(\H)}}

\newcommand{\ca}{\mathrm{C}^*}
\newcommand{\cstar}{\mathrm{C}^*}
\newcommand{\cenv}{\mathrm{C}^*_{\text{env}}}

\newcommand{\bbC}{{\mathbb{C}}}

\newcommand{\bbN}{{\mathbb{N}}}

  \newcommand{\A}{{\mathcal{A}}}
  \newcommand{\B}{{\mathcal{B}}}

\renewcommand{\H}{{\mathcal{H}}}
  \newcommand{\I}{{\mathcal{I}}}
  \newcommand{\J}{{\mathcal{J}}}

  \newcommand{\N}{{\mathcal{N}}}

\renewcommand{\S}{{\mathcal{S}}}
  \newcommand{\T}{{\mathcal{T}}}
  \newcommand{\U}{{\mathcal{U}}}

  \newcommand{\X}{{\mathcal{X}}}

\renewcommand{\phi}{\varphi}

\newcommand{\fA}{{\mathfrak{A}}}

\newcommand{\fB}{{\mathfrak{B}}}

\newcommand{\fI}{{\mathfrak{I}}}

\newcommand{\fM}{{\mathfrak{M}}}

\newcommand{\Alg}{\operatorname{Alg}}

\newcommand{\dirlim}{\varinjlim}

\newcommand{\Id}{\operatorname{Id}}
\newcommand{\LL}{\operatorname{L}}

\newcommand{\ran}{\operatorname{Ran}}
\newcommand{\rank}{\operatorname{rank}}

\newcommand{\sumoplus}{\operatornamewithlimits{\sum\oplus}}
\newcommand{\ba}{\backslash}
\begin{document}

%%%%%%%%%%%%%%%%%%%%%%%%%%%%%%%%%%%%%%%%%%%%%%%%%%%%%%%%%%%%
\title[Meet irreducible ideals]{Meet irreducible ideals and
representations of limit algebras}
\thanks{October 20, 2001 version.}
\author[K.R.Davidson]{Kenneth R. Davidson}
\thanks{First author partially supported by an NSERC grant.}
\address{Pure Math.\ Dept.\\U. Waterloo\\Waterloo, ON\;
N2L--3G1\\CANADA}
\email{krdavids@uwaterloo.ca}
\author[E. Katsoulis]{Elias~Katsoulis}
\thanks{Second author's research partially supported by a grant from
ECU}
\address{Math.\ Dept.\\East Carolina University\\
Greenville, NC 27858\\USA}
\email{KatsoulisE@mail.ecu.edu}
\author[J. Peters]{Justin~Peters}
\address{Math.\ Dept.\\Iowa State University\\Ames, IA
50011\\USA}
\email{peters@iastate.edu}
\begin{abstract}
In this paper we give
criteria for an ideal $\J$ of a TAF algebra $\A$ to be meet irreducible.
We show that $\J$ is meet irreducible if and only if the
$\ca$-envelope
of $\A/\J$ is primitive. In that case, $\A/\J$ admits a faithful nest
representation which extends to a $*$-representation of the
$\ca$-envelope for $\A / \J $. We
also characterize the meet irreducible ideals as the kernels of 
bounded nest representations;
this settles the question of whether the n-primitive and
meet-irreducible ideals coincide.
\end{abstract}

\subjclass{47L80}
\date{}
\maketitle
%%%%%%%%%%%%%%%%%%%%%%%%%%%%%%%%%%%%%%%%%%%%%%%%%%%%%%%%%%%%
\section{Introduction} \label{s:intro}
Representation theory of operator algebras is still in its infancy.
While for $\cstar$-algebras the fundamentals of
representation theory have long been known, for nonselfadjoint algebras
there are hardly any results of a general nature.
For `triangular operator algebras' (a term which we leave undefined),
intuition suggests that the fundamental building
blocks for representation theory should be nest representations. In the
category of $\cstar$-algebras,
the nest representations are precisely the irreducible representations.

Recall that a \emph{nest representation} is a representation for which
the closed, 
invariant subspaces form a nest (i.e., are linearly ordered).
In his study of nonselfadjoint crossed products, Lamoureux introduced
the notion of \emph{n-primitive ideal}. An
ideal is n-primitive if it is the kernel of a nest representation.
Lamoureux has shown that in various contexts
in nonselfadjoint algebras the n-primitive ideals play a role analogous
to the primitive ideals in $\cstar$-algebras.
Thus, one can give the set of n-primitive ideals the hull-kernel
topology, and for every (closed, two-sided) ideal
$\I$ in the algebra, $\I$ is the intersection of all n-primitive ideals
containing $\I;$
in other words,  $ \I = k(h(\I)).$

An ideal $\J$ of an algebra $\A$ is \textit{meet irreducible} if, for
any
ideals $\I_1$ and $\I_2$ containing $\J,$ the relation
 $\I_1 \cap \I_2 = \J$ implies that either $\I_1
= \J$ or $\I_2 = \J$. In the case of $T_n,$ the algebra of upper
triangular $n\times n$ matrices,
meet irreducible ideals are gotten by `cutting a wedge' from the
algebra:  let $ 1 \leq i_0 \leq j_0 \leq n.$ The ideal
\[ \I = \{ (a_{ij}):  a_{ij} = 0,\ i_0 \leq i \leq j \leq j_0 \} \]
is meet irreducible, and every meet irreducible ideal of $T_n$ has this
form.

The relationship between
meet-irreducible and n-primitive ideals is studied in
a variety of examples in \cite{La97}, and in \cite{DHHLS}
meet-irreducible ideals in strongly maximal
triangular AF-algebras are characterized by sequences of matrix units
and also in terms of groupoids.  In that paper
it is shown that every meet-irreducible ideal is n-primitive. This is
done by constructing a nest representation.
The converse
question, whether every
n-primitive ideal is meet-irreducible, was left open.

In a recent work \cite{DK}, the first
two authors examined the $\cstar$-envelope of a quotient $\A/\J$
of a strongly maximal TAF algebra by an ideal $\J$.  They showed
the $\cstar$-envelope is an AF $\cstar$ algebra, even though the
quotient $\A/\J$ is not in general a TAF algebra.
It turns out that the $\cstar$-envelope of $\A/\J$ is sensitive enough
to detect the meet irreducibility of $\J$. In
Theorem~\ref{mi characterization} we show that $\J$ is meet irreducible
if and only if
$\cenv (\A/\J)$, the $\cstar$-envelope of $\A/\J$, is primitive. The
theory of $\ca$-envelopes provides the natural framework
for studying results of this type. In Theorem~\ref{mi repn} we show that
for a meet irreducible ideal $\J$,
there exists a faithful and irreducible $*$-representation of $\cenv
(\A/\J)$, whose restriction on $\A / \J$ is a
nest representation. Since the converse is easily seen to be true,
Theorem~\ref{mi repn}
provides a characterization of meet irreducible ideals in terms of the
representation theory for $\A$.

The question of whether the kernel of a
nest representation is a meet-irreducible ideal emerged at the
Ambelside, U.K.
conference in summer, 1997. Subsequently some progress was made. In
\cite{HPP00} a partial result was
obtained: if the TAF algebra $\A$ has totally ordered spectrum, or if
the
nest representation $\pi$ has the property that
the von Neumann algebra generated by $\pi(\A \cap \A^*)$ contains an
atom,
then ker$(\pi)$ is meet-irreducible. The solution presented in Theorem
\ref{t:nestrep} is self-contained and does not make use of the results
of
\cite{DHHLS} or \cite{HPP00}.  Thus the question is now settled for
strongly maximal TAF algebras.

Despite the fact that evidence at hand is limited, it nonetheless 
seems worthwhile to ask
\begin{ques}
Are there any operator algebras for which the n-primitive ideals,
and the meet-irreducible ideals do not coincide?
\end{ques}

We would like to thank Alan Hopenwasser for several suggestions,
including one leading to the final version of
Theorem \ref{t:nestrep}.

\section{The main results}

We begin by recalling a result of Lamoureux \cite{La97}.

\begin{lem} \label{Lam}
Let $\I$ be a closed, two-sided ideal in a separable $\cstar$-algebra
$\A.$
Then the following are equivalent:
\begin{enumerate}
\item $\I$ is n-primitive
\item $\I$ is primitive
\item $\I$ is prime
\item $\I$ is meet-irreducible
\end{enumerate}
\end{lem}

One can actually characterize when an AF $\ca$-algebra
is primitive in terms of its Bratteli diagram.
Let $\fA = \dirlim (\fA_i, \phi_i)$ be an AF $\ca$-algebra
and assume that each $\fA_i$ decomposes as a direct sum
$\fA_i = \oplus_{j}\fA_{i\, j}$ of finite dimensional
full matrix algebras $\fA_{i \, j}$. A \textit{path} $\Gamma$
for $\fA = \dirlim (\fA_i, \phi_i)$ is a sequence 
$\{ \fA_{i \, j_{i}} \}_{i=1}^{\infty}$
so that for each pair of nodes
$\left( (i , j_i), (i+1, j_{i+1}) \right)$ there exist
an arrow in the Bratteli diagram for $\fA = \dirlim (\fA_i, \phi_i)$
which joins them. It is known that $\fA$ is primitive iff there is a
path
$\Gamma$
for $\fA = \dirlim (\fA_i, \phi_i)$ so that each summand of $\fA_i$ is
eventually
mapped into a member of $\Gamma$. We call such a path $\Gamma$ an
\textit{essential
path} for $\fA$.

Beyond $\ca$-algebras, a meet irreducible ideal need not be primitive. 
In \cite{DHHLS}, a description of all meet irreducible ideals of a
TAF algebra was given in terms of matrix unit sequences.

\begin{defn} \label{mi-chain}
Let $\A = \dirlim (\A_i, \phi_i)$ be a TAF algebra. A sequence $(
e_i)_{i
\geq N}$ of matrix units
from $\A$ will be called an \textit{mi-chain} if the following two
conditions are satisfied for all
$i \geq N$:
\begin{itemize}
\item[ (A)] $e_i \in \A_i$.
\item[ (B)] $e_{i+1} \in \Id_{\, i+1}(e_i)$,
\end{itemize}
where $\Id_{\, i+1}(e_i)$ denotes the ideal generated by $e_i$ in
$\A_{i+1}$.
\end{defn}

If $( e_i)_{i \geq N}$ is an mi-chain for $\A = \dirlim (\A_i, \phi_i)$,
let $\J$ be the join of all ideals which do not
contain any matrix unit $e_i$ from the chain.
In \cite[Theorem 1.2]{DHHLS} it is shown that for a TAF algebra $\A =
\dirlim (\A_i, \phi_i)$,
given an mi-chain $( e_i)_{i \geq N}$, the ideal $\J$ associated with $(
e_i)_{i \geq N}$ is meet
irreducible. Conversely, every proper meet irreducible ideal in $\A =
\dirlim
(\A_i, \phi_i)$
is induced by some mi-chain, chosen from some contraction of this
representation.

In this paper we give a characterization of the meet irreducible ideals
of
TAF algebras
in terms of $\ca$-envelopes of quotient algebras. We need to recall the
notation and
machinery from \cite{DK}.

Let $\fA = \dirlim (\fA_i, \phi_i)$ be the enveloping $\ca$-algebra for
a TAF algebra $\A = \dirlim (\A_i, \phi_i)$.  Let $\J \subseteq \A$
be a closed ideal, and let $\J_i := \J \cap \A_i$.
For each $i\ge1$, $\S_i$ denotes the collection of all diagonal
projections
$p$
which are semi-invariant for $\A_i$, are supported on a single
summand of $\fA_i$ and satisfy $(p \A_i p) \cap \J = \{ 0 \}$. We form
finite dimensional C*-algebras
\[
 \fB_i := \sumoplus_{p \in S_i} \B(\ran p)
\]
where $\B(\ran p)$ denotes the bounded operators on $\ran p$.
Of course, $\B(\ran p)$ is isomorphic to $\fM_{\rank p}$.
Let $\sigma_i$ be the map from $\fA_i$ into $\fB_i$
given by $\sigma_i(a) = \sum^\oplus_{p \in S_i} pap|_{\ran p}$.
The map $\sigma_i|_{\A_i}$ factors as
$\rho_i q_i$ where $q_i$ is the quotient map of $\A_i$ onto
$\A_i/ \J_i$ and $\rho_i$ is a completely isometric homomorphism of
$\A_i/ \J_i$ into $\fB_i$. 
Notice that $\fB_i$ equals the $\ca$-algebra generated by $\rho_i ( \A_i/
\J_i )$.

We then consider unital embeddings $\pi_i$ of $\fB_i$ into $\fB_{i+1}$
defined as follows.
For each $q \in \S_{i+1}$ we choose projections $p \in \S_i$ which
maximally
embed into $q$ under the action of $\phi_i$. This way, we determine
multiplicity one
embeddings of $\B (\ran p)$ into $\B (\ran q)$. Taking into account all
such
possible embeddings,
we determine the embedding $\pi_i$ of $\fB_i$ into $\fB_{i+1}$.

Finally we form the subsystem of the directed limit $\fB= \dirlim
(\fB_i,
\pi_i)$
corresponding to all summands which are \textit{never} mapped into a
summand $\B(\ran p)$ where $p$ is a maximal element of some $\S_i$.
Evidently this system is directed upwards.
It is also hereditary in the sense that if every image of a summand
lies in one of the selected blocks, then it clearly does not map into
a maximal summand and thus already lies in our system.
By \cite[Theorem~III.4.2]{D}, this system
determines an ideal $\fI$ of $\fB$.
The quotient $\fB' = \fB/ \fI$ is the AF algebra corresponding to the
remaining summands and the remaining embeddings; it can be expressed as
a
direct
limit $\fB^{\prime} = \dirlim (\fB_{i}^{\prime}, \pi_{i}^{\prime})$,
with
the understanding that
$\fB_{i}^{\prime}= \oplus_{j} \fB_{i\, j}$ for these remaining summands
$\fB_{i\, j}$
of $\fB_i$. It can be seen
that the quotient map is isometric on $\A / \J$ and that $\fB^\prime$
is the $\ca$-envelope of $A/ \J$.

\begin{thm} \label{mi characterization}
Let $\A$ be a TAF algebra and let $\J \subseteq \A$ be an ideal. Then
$\J$
is meet
irreducible if and only if the algebra $\cenv (\A/\J)$ is primitive.
\end{thm}

\Prf Assume that $\fB^{\prime} = $ $\cenv (\A/\J)$ is primitive and let
$\Gamma=\left( \fB_{i \, j_{i}} \right) _{i=1}^{\infty}$ essential path
for $\fB^{\prime}$. Let $e_i$ for $\fB_{i \, j_{i}}$ be the
\textit{characteristic matrix units}
for $\fB_{i \, j_{i}}$, i.e., the ones on the top right corner of
$\fB_{i \,
j_{i}}$.

Assume that there exist ideals
$\I_1$ and $\I_2$, properly containing $\J$.
Since
$\I_1$ and $\I_2$  properly contain $\J$, there
exist matrix units $f_k \in \I_k$ with $f_k \notin \J$, $k =1,2$. So the
images
of the $f_k$ appear in the presentation for the $\ca$-envelope in
perhaps different summands. However, the existence of an essential path
$\Gamma$
implies that some
subordinates for the $f_k$ will appear in some member of $\Gamma$, say
$\fB_{i \, j_{i}}$, for $i$ sufficiently large,
and so
$e_i \in  \I_1 \cap \I_2$. However, $e_i \notin \J$ and so $\J$ is
properly
contained in
$ \I_1 \cap \I_2$. It follows $\J$ is meet irreducible.

Conversely, assume that $\J$ is meet irreducible. In light of
Lemma~\ref{Lam} and the
subsequent comments, it suffices to show that the trivial ideal $\{ 0
\}$ is
meet
irreducible in the $\ca$-envelope $\fB^{\prime}$.

By way of contradiction assume that there are non-trivial ideals
$\I_1$ and $\I_2$ of $\fB^{\prime}$ so that $ \I_1 \cap \I_2 = \{ 0 \}$.
We
claim that
$(\A / \J) \cap \I _k \neq \{ 0 \}$, $k = 1, 2$. Indeed, any non-trivial
summand of  $\I_k$ will eventually be mapped into a direct summand
$\fB_{i \, j_{i}}$ of $\fB^{\prime}$ corresponding to some maximal
element
of $\S_i$. Hence all matrix units in $\fB_{i \, j_{i}}$ belong to
$\I_k$,
including the characteristic one. This one however also belongs to $\A /
\J$
and therefore
in the intersection $(\A / \J) \cap \I _k$.

The claim shows that the
zero ideal is not meet irreducible in $\A / \J$. By considering the
pullback, this implies that $\J$ is not meet irreducible in $\A$,
which is the desired contradiction.
\bx

Notice that the sequence $( e_i)_{i=1}^{\infty}$ associated with the
path
$\Gamma$
in the proof above satisfies the Conditions (A) and (B) of the
Definition~\ref{mi-chain}
and is therefore an mi-chain for the ideal $\J$.

\begin{thm} \label{mi repn}
If $\A$ is a TAF algebra and $\J$ an ideal of $\A$, then the following
are
equivalent:
\begin{itemize}
\item[(i)] There exists a faithful representation $\tau : \cenv (\A /
\J)
\longrightarrow \BH $ so that
$\tau (\A / \J)$ is weakly dense in some nest algebra.
\item[(ii)] $\J$ is meet irreducible.
\end{itemize}
\end{thm}

\Prf Assume that (i) is valid and let $\tau : \cenv (\A / \J)
\longrightarrow \BH $ be a faithful
representation so that
$\tau (\A / \J)$ is weakly dense in some nest algebra $\Alg \N$. By way
of
contradiction assume that
$\J$ is not meet irreducible. Theorem~\ref{mi characterization} and
Lemma~\ref{Lam}
imply the existence of nonzero closed ideals
$\I_1$ and $\I_2$ in $ \cenv (\A / \J)$ so that $\I_1 \I_2 = \{ 0 \}$.
The non-zero subspaces $[\tau (\I_i)\H]$ must be mutually
orthogonal. However they are both invariant under $\tau (\A / \J)$,
and therefore belong to $\N$, a contradiction.

Conversely, assume that (ii) is valid and so, by Theorem~\ref{mi
characterization}, we know that
$ \cenv (\A / \J)$ is primitive.
Retain the notation established in the paragraphs
preceding Theorem~\ref{mi characterization}. Hence
\[
  \cenv (\A / \J)=  \fB^{\prime} = \dirlim (\fB_{i}^{\prime},
\pi_{i}^{\prime})
\]
where $\fB_{i}^{\prime}= \oplus_{j} \fB_{i\, j}$ for the remaining
summands
$\fB_{i\, j}$
of $\fB_i$. Let $\Gamma=\left( \fB_{i \, j_{i}} \right) _{i=1}^{\infty}$
be the essential path
for $\fB^{\prime}$

Each $\fB_{i\, j}$ is a full matrix algebra and therefore contains the
algebra $\T_{i \, j}$ of upper triangular matrices. Form the finite
dimensional algebras $\T_{i}^{\prime}= \oplus_{j} \T_{i\, j}$ and
consider the direct limit algebra
\[
\T^{\prime} = \dirlim ( \T^{\prime}_{i}, \pi_{i}^{\prime}),
\]
where $\pi_{i}^{\prime}$ is as earlier. Clearly, $\T^{\prime}$ is a TAF
algebra whose
enveloping $\ca$-algebra is $\fB^{\prime}$. Moreover, $\T^{\prime}$
contains
$\A /\J$.

We define a state $\omega$ on $\fB^{\prime}$ as follows. Let
$(p_i)_{i=1}^{\infty}$ be a sequence
of diagonal projections with $p_i \in \T_{i \, j_{i}}$ so that $p_{i+1}$
is
a subordinate of $p_i$,
$i \in \bbN$. We define $\omega _i : \fB_i^{\prime} \rightarrow \bbC$ to
be
the compression on $p_i$
and we let
$\omega$ to be the direct limit  $\omega = \dirlim \omega _i$. Consider
the
GNS triple $( \tau , \H, g)$
associated with the state $\omega$, i.e., $\tau$ is a representation of
$\fB^{\prime}$ on $\H$ and
$g \in \H$ so that $\omega (a) =\langle \tau ( a ) g, g\rangle, a \in
\fB^{\prime}$. Since $\omega$ is
pure, $\tau$ is irreducible. Moreover, $p_i \in \fB_{i \, j_{i}}, i \in
\bbN$ and so $\tau$ is also faithful.

An alternative
presentation for $( \tau, \H, g)$ was given in \cite[Proposition
II.2.2]{OPe}. Since $\omega$ is multiplicative on
the diagonal $\T^{\prime} \cap (\T^{\prime})^{*}$,
one considers  $\H$ to be $\LL^{2} (\X ,\mu)$,
where $\X$ is the Gelfand spectrum of $\T^{\prime} \cap
(\T^{\prime})^{*}$
and $\mu$ the counting measure on the orbit of $\omega$ in $\X$. With
these
identifications, given any matrix unit $e$,
$\tau ( e)$ is the translation operator on $\X$ defined in the
paragraphs
preceding \cite[Theorem II.1.1]{OPe}.

In \cite[Proposition II.2.2]{OPe} it is shown that $\tau$ maps
$\T^{\prime}$
onto a weakly dense subset of some nest algebra. 
The proof of the theorem will follow if we show that the weak closure
of $\tau ( \A / \J)$ contains  $\tau (\T^{\prime})$.

A moment's reflection shows that given any contraction 
$a \in \T^{\prime}_{i\, j}$
and matrix units $e_1, e_2, \dots, e_n $ and $f_1, f_2, \dots, f_n$
in $ \fB^{\prime}$, there exists a contraction
$\hat{a} \in \A / \J$ so that
\[
 \omega (f_{k}^{*} \hat{a} e_{k} )= \omega (f_{k}^{*} a e_{k})
\]
and therefore
\[
 \langle \tau (\hat{a})\tau(e_{k})g, \tau(f_{k})g \rangle =
 \langle  \tau (a) \tau(e_{k})g ,\tau( f_{k})g \rangle
\]
for all $k = 1, 2, \ldots, n$. However the collection of all vectors of
the
form $\tau(e)g$, where $e$ ranges over all matrix units of
$\fB^{\prime}$, forms a dense subset of $\H$ and
so the desired density follows.
\bx

\begin{rems} \ 
\begin{enumerate}
\renewcommand{\labelenumi}{(\arabic{enumi}) }
\item The implication (i) $\implies$ (ii) also follows from
Theorem \ref{t:nestrep}. 
\item There exists a faithful representation $\tau$ of $\cenv(\A/\J)$
in $\B(\H)$ so that $\tau(\A/\J)$ is weakly dense in a nest algebra if
and only if there is a faithful \textit{irreducible} representation
$\tau$ of $\cenv(\A/\J)$ in $\B(\H)$ so that $\tau(\A/\J)$ is weakly
dense in a nest algebra.
\end{enumerate}
\end{rems}

\begin{thm} \label{t:nestrep}
Let $\A$ be a strongly maximal TAF algebra, and let $\pi$ be a bounded
nest representation of $\A$ on a Hilbert space $\H$.  
Then ker$(\pi)$ is a meet-irreducible ideal.
\end{thm}

\Prf
Since a bounded representation of the diagonal masa $\A \cap \A^* $ is 
completely bounded, c.f.\ \cite[Theorem 8.7]{Paul}, and
a completely bounded representation is similar to a completely
contractive representation by \cite[Theorem 8.1]{Paul}, 
we may assume that the restriction of  $\pi$ to the diagonal masa is
completely contractive. It follows
that the restriction of $\pi$ to the diagonal masa is a
star representation.
Let $\J = \text{ker}(\pi),$ and $\J_1, \J_2$ be ideals in $\A$ properly
containing $\J$.
We need to show that $\J_1 \cap \J_2$ properly contains $\J$.

Since $\pi$ is a nest representation, we have (after possibly 
interchanging $\J_1$ and $J_2$) that
\[
(0) \neq [\pi(\J_1)\H] \subseteq [\pi(\J_2)\H]
\]
where $[\X]$ denotes the closed subspace generated by $\X \subset \H$. 
Fix $n \in \bbN$ and let $u$ be a matrix unit in $ \J_1 \cap \A_n \ba \J $.
Choose a vector $h \in \H$ be such that $ \|\pi(u)h \| = 1$. 
There exist $m \geq n$, $N \ge 1$, matrix units $v_t \in \J_2\cap \A_m$
and vectors $h_t \in \H$ for $1 \leq t \leq N$ such that
\[
  \Big\| \pi(u)h - \sum_{t=1}^N \pi(v_t) h_t \Big\| < \frac{1}{4} \ .
\]
In particular, $\big\| \sum_{t=1}^N \pi(v_t) h_t \big\| > 3/4$.
We may assume that $\pi(v_t) \neq 0$ for all $t$.

Define a projection $E = \bigvee_{t=1}^N \pi(e_t) =
\pi \big( \bigvee_{t=1}^N e_t \big)$, where $e_t = v_tv_t^*$ are
diagonal matrix units (which need not be distinct). 
For all $s, t,$
\[
 \pi(e_s)\pi(v_t) = \pi(v_sv_s^*v_t) =
 \begin{cases}
 \pi(v_t) &\text{ if } v_sv_s^* = v_tv_t^* \\
        0 &\text{ otherwise }
 \end{cases}
\]
So we have $E\sum\pi(v_t)h_t = \sum \pi(v_t)h_t $.
Now
\begin{align*}
  \Big\| E\pi(u)h - \sum_{t=1}^N \pi(v_t)h_t \Big\| &=
  \Big\| E \big( \pi(u)h - \sum_{t=1}^N \pi(v_t)h_t \big) \Big\|
 \\ &\le \Big\| \pi(u)h - \sum_{t=1}^N \pi(v_t)h_t \Big\| < \frac14 .
\end{align*}
Therefore $\| E \pi(u)h \| > 1/2$.
In particular, there exists at least one $t$, $ 1 \leq t \leq N$, such
that $\pi(e_t)\pi(u) \neq 0$.

Embed $u \in \A_n \hookrightarrow \A_m$ and decompose it as a sum $u =
\sum u_s$ of matrix units
in $\A _m$. Then
$e_t u = u_s$, for some $s$, and so it follows $\pi(u_s) \neq
0$, i.e., $u_s \notin \J$.
Thus we have identified matrix units $u_s \in \J_1 \ba \J$ and $ v_t \in
\J_2 \ba \J$ of $\A_m$ with the same final
projection. Say $u_s = e^{(m,r)}_{ij}$ and $v_t = e^{(m,r)}_{ik}$. We
now distinguish three
cases:
\begin{equation*}
   \begin{aligned}
    \text{If } j &= k, \text{ then } u_s = v_t \in \J_1 \cap \J_2 \ba \J
;\\
    \text{If } j &< k, \text{ then }  v_t = u_s e^{(m,r)}_{jk} \in \J_1
\cap \J_2 \ba \J;\\
    \text{If } j &> k, \text{ then }  u_s = v_t e^{(m,r)}_{kj} \in \J_1
\cap \J_2 \ba \J.
   \end{aligned}
\end{equation*}
It follows that in all three cases $\J_1 \cap \J_2$ properly contains
$\J$.
Thus $\J$ is meet-irreducible.
\bx

%%%%%%%%%%%%%%%%%%%%% References %%%%%%%%%%%%%%%%%%%%%


\begin{thebibliography}{99}

\bibitem{D} K. Davidson,
\textit{$\ca$-algebras by example},
Fields Institute Monographs, American Mathematical Society, 1996.

\bibitem{DK} K. Davidson and E. Katsoulis,
\textit{Primitive limit algebras and $\ca$-envelopes},
Adv.\ Math., \textit{to appear}

\bibitem{DHHLS} A. Donsig, A. Hopenwasser, T. Hudson, M. Lamoureux
and B. Solel, \textit{Meet irreducible ideals in direct limit algebras}
Math. Scand., \textbf{87}, 2000, 27--63.

\bibitem{HPP00}Alan Hopenwasser, Justin Peters, Stephen Power,
\textit{Nest Representations of TAF Algebras},
Canad. J. Math. \textbf{52}, 2000, 1221--1234.
 
\bibitem{La93} Michael Lamoureux,
\textit{Nest Representations and Dynamical Systems},
J. Funct. Anal. \textbf{114}, 1993, 345--376.

\bibitem{La96} Michael Lamoureux,
\textit{Ideals in some continuous nonselfadjoint crossed product
algebras},
J. Funct. Anal. \textbf{142}, 1996, 211--248.

\bibitem{La97} Michael Lamoureux,
\textit{Some Triangular AF Algebras},
J. Operator Theory \textbf{37}, 1997, 91--109.

\bibitem{OPe}  J. Orr and J. Peters,
\textit{Some representations of TAF algebras},
Pacific.\ J.\ Math. \textbf{167} (1995), 129--161.

\bibitem{Paul} V. I. Paulsen,
\textit{Completely Bounded Maps and Dilations},
Longman Scientific, New York, Wiley, 1986.

\bibitem{Ped} G. K. Pedersen,
\textit{C$^*$-Algebras and their automorphism groups},
London Mathematical Society Monograph, Academic Press, 1979.

\bibitem{SV} \c S. Str\u atil\u a and D.V. Voiculescu,
\textit{Representations of AF-algebras and of the group $U(\infty)$},
Springer Lect.\ Notes Math.\ \textbf{486},
Springer Verlag, Berlin, New York, 1975.

\bibitem{Thel91} Michael Thelwall,
\textit{Dilation theory for subalgebras of AF algebras},
J. Operator Theory \textbf{25}, 1991, 275--282.

\end{thebibliography}
\end{document}